\documentclass[12pt]{amsart}
\usepackage{amsmath,amssymb,url}
\usepackage{hyperref}
\usepackage{amsthm}
\usepackage{graphicx}
\usepackage{tikz}
\usepackage{tikz-cd}
\usepackage{stmaryrd}
\usepackage{mathrsfs,upgreek}
\usepackage{eucal}
\usepackage{enumerate}
\usepackage{changes}
\usepackage{verbatim}

\newtheorem{Thm}{Theorem}[section] 
\newtheorem{Alg}[Thm]{Algorithm}
\newtheorem{Prop}[Thm]{Proposition} 
\newtheorem{Lem}[Thm]{Lemma} 
 
\theoremstyle{definition}
\newtheorem{Rem}[Thm]{Remark} 
 
\theoremstyle{definition}
\newtheorem{Def}[Thm]{Definition}
\numberwithin{equation}{section}

\newcommand{\join}{\vee}

\newcommand{\alg}[1]{{\textbf{\upshape #1}}}  %
\newcommand{\vv}[1]{\mathcal {#1}}
\renewcommand{\a}{\alpha}
\renewcommand{\b}{\beta}

\newcommand{\f}{\varphi}
\newcommand{\g}{\gamma}

\renewcommand{\t}{\tau}

\newcommand{\con}{\operatorname{Con}}

\newcommand{\Pol}{\operatorname{Pol}}

\newcommand{\Eqv}{\mathrm{Eqv}}

\newcommand{\Eq}{\mathrm{Eq}}
\newcommand{\Cg}{\mathrm{Cg}}

\newcommand{\NN}{{\mathbb{N}}}

\begin{document}
	\markboth{Stefano Fioravanti}
	{Commutator equations}
	
	\title[Commutator equations]{On properties described by terms in commutator relation}

	\author{Stefano Fioravanti}
	
	\address{Stefano Fioravanti,
		DIISM, 
		Universit\`a di Siena,
		Via Roma 56, 
		53100 Siena,
		Italy}
	\email{\tt stefano.fioravanti66@gmail.com}
	
	\begin{abstract}
		We investigate properties of varieties of algebras described by a novel concept of equation that we call \emph{commutator equation}. A commutator equation is a relaxation of the standard term equality obtained substituting the equality relation with the commutator relation. Namely, an algebra $\alg A$ satisfies the commutator equation $p \approx_{C} q$ if for each congruence $\theta$ in $\con (\alg A)$ and for each substitution $p^{\alg A}, q^{\alg A}$ of elements in the same $\theta$-class, then $(p^{\alg A}, q^{\alg A}) \in [\theta, \theta]$. This notion of equation draws inspiration from the definition of \emph{weak difference term} and allows for further generalization of  it.
		Furthermore, we present an algorithm that establishes a connection between congruence equations valid within the variety generated by the abelian algebras of the idempotent reduct of a given variety and congruence equations that hold within the entire variety.
		 Additionally, we provide a proof that if the variety generated by the abelian algebras of the idempotent reduct of a variety satisfies a non-trivial idempotent Mal'cev condition then also the entire variety satisfies a non-trivial idempotent Mal'cev condition, statement that follows also form \cite[Theorem 3.13]{KK.TSOC}. 
	\end{abstract}

\maketitle
	
	\section{Introduction}
	
	The genesis of this work is our interest in properties of classes of varieties described by equations. Starting from Mal'cev's influential work on congruence permutability as presented in \cite{Mal.OTGT}, the problem of characterizing properties of classes of varieties using Mal'cev conditions has yielded numerous results. Notably, A. Pixley in \cite{Pix.DAPO} identified a strong Mal'cev condition that describes the class of varieties with distributive and permuting congruences. Similarly, B. Jónsson in \cite{Jon.AWCL} established a Mal'cev condition that characterizes congruence distributivity, and A. Day in \cite{Day.ACOM} provided a Mal'cev condition that characterizes the class of varieties with modular congruence lattices.
	
	The aforementioned results are examples of a more general theorem obtained independently by Pixley \cite{Pix.LMC} and R. Wille \cite{Wil.K} that can be considered as a foundational result in the field. They proved that if $p \leq q$ is a lattice identity, then the class of varieties whose congruence lattices satisfy $p \leq q$ is the union of countably many Mal'cev classes.  In particular, \cite{Pix.LMC} and \cite{Wil.K} provide an algorithm for generating Mal'cev conditions associated with congruence identities.
	
	Furthermore, these investigations have led to the study of equations where variables range not only over congruence lattices but also over potentially broader sets, such as the lattices of all tolerances or all compatible reflexive relations. An analogous of  Pixley \cite{Pix.LMC} and R. Wille \cite{Wil.K} result for equations over the lattice of compatible reflexive relations of an algebra can be found in \cite{Fio.MCCT}.
	
	Moreover, the study of Mal'cev conditions has gained significant popularity in recent years due to its profound association with Constraint Satisfaction Problems (CSPs) \cite{Bul.ADTF,Zhuk.APOT} and is closely connected also to the study of properties of closed sets of operations, (e.g. \cite{Fio.CSOF2, Fio.EOAS, Spa.OTNC}).
	
	Following this branch of research, our aim is to introduce a new type of equations, that we will call \emph{commutator equations}. Drawing inspiration from the definition and characterization of varieties with a \emph{weak difference} term (as discussed in Chapter $6.1$ of \cite{KK.TSOC}), we formulate the concept of commutator equations to generalize this notion. It is evident that the \emph{weak difference} term is a relaxation of a Mal'cev term obtained by loosening the characterizing equations to what we refer to as commutator equations. Consequently, this produces a property defining a Mal'cev condition that is strictly weaker than the one characterizing the class of varieties with a Mal'cev term. 
	
	In Section \ref{Notations} and \ref{LabelledG} we introduce some of the basic concepts for our investigation. In Section \ref{Taylor} we prove an equivalent of the Taylor's theorem \cite{Tay.VOHL} for commutator equations, result that also follows from \cite[Theorem 3.13]{KK.TSOC}. Namely,  Theorem \ref{TaylorThm} shows that a variety $\vv{V}$ is Taylor if and only if  the variety $\vv{V}'$ generated by the abelian algebras of the idempotent reduct of $\vv{V}$ is Taylor. Furthermore, we provide a new characterization of this property in terms of  congruence equations. This result establishes a connection between a variety $\vv V$ and the variety $\vv{V}'$ generated by the abelian algebras of the idempotent reduct of $\vv{V}$ that we develop further in Section \ref{Alg}. In this last section our aim is to provide a Pixley-Wille type of algorithm for commutator equation. In the main result of the section, Theorem \ref{ThmweakAlg}, we prove that, under mild assumptions on the lattice terms involved, weakening standard term equations  produced by a congruence equation via the Pixley-Wille algorithm to commutator equations  gives  a property describing a weak Mal'cev class. 
	
	\section{Preliminaries and Notations}{}
	\label{Notations}
	
	In this section, we provide a brief review of the fundamental definitions in universal algebra, along with the introduction of a few new definitions that are relevant to our study. For basic concepts in general algebra, such as lattices, algebras, and varieties, we refer to the textbook \cite{BS.ACIU}. For more advanced topics, including abelian congruences and the commutator of congruences, we direct the reader to \cite{KK.TSOC}. The general theory of Mal'cev conditions and Mal'cev classes can be found in the classical treatment presented in \cite{Taylor1973}, or alternatively, in the more modern approach outlined in \cite{KK.TSOC}.
	We begin by recalling the definition of  \emph{Centralizer}.  Let $\mathbf{A}$ be an algebra. For $\alpha, \beta, \delta  \in \con(\mathbf{A})$ we say that $\alpha$ \emph{centralizes} $\beta$ \emph{modulo} $\delta$ if for all possible matrices
    \begin{center}
		\begin{equation*}
			\begin{pmatrix}
				f(\mathbf{a},\mathbf{u}) & f(\mathbf{a},\mathbf{v}) \\ f(\mathbf{b},\mathbf{u}) & f(\mathbf{b},\mathbf{v})
			\end{pmatrix}
		\end{equation*}
	\end{center}
 %For $\alpha, \beta,  \in \con(\mathbf{A})$ we define $M(\alpha,\beta)$ as the set of matrices of the form 
%	\begin{center}
%		\begin{equation*}
%			\begin{pmatrix}
%				f(\mathbf{a},\mathbf{u}) & f(\mathbf{a},\mathbf{v}) \\ f(\mathbf{b},\mathbf{u}) & f(\mathbf{b},\mathbf{v})
%			\end{pmatrix}
%		\end{equation*}
%	\end{center}	
	with $f \in \Pol_{n+m}(\mathbf{A})$, $\mathbf{a}, \mathbf{b} \in A^n$, $\mathbf{u}, \mathbf{v} \in A^m$ with $\mathbf{a} \equiv_{\alpha} \mathbf{b}$ and $\mathbf{u} \equiv_{\beta} \mathbf{v}$,	
	%For $\alpha, \beta, \delta \in \con(\mathbf{A})$, we say that $\alpha$ \emph{centralizes} $\beta$ \emph{modulo} $\delta$ if for all possible matrices $M \in M(\alpha,\beta)$ with $f \in \Pol(\mathbf{A})$ 
 we have that:
	\begin{center}
		if $f(\mathbf{a},\mathbf{u}) \equiv_{\delta} f(\mathbf{a},\mathbf{v})$, then $f(\mathbf{b},\mathbf{u}) \equiv_{\delta} f(\mathbf{b},\mathbf{v})$.
	\end{center}
	In this case we write \index{$C(\alpha,\beta;\delta)$}$C(\alpha,\beta;\delta)$ and we call the ternary relation $C$ \emph{centralizer}. For a comprehensive overview of the various properties of the centralizer relation, we refer the interested reader to \cite[Section $2.5$]{KK.TSOC}.
	
	Let us now recall the definition of the \emph{two-term Centralizer}. Let $\mathbf{A}$ be an algebra. For $\alpha, \beta, \delta  \in \con(\mathbf{A})$ we say that $\alpha$ \emph{two-term centralizes} $\beta$ \emph{modulo} $\delta$ if for all possible matrices
		\begin{center}
		\begin{equation*}
			\begin{pmatrix}
				f(\mathbf{a},\mathbf{u}) & f(\mathbf{a},\mathbf{v}) \\ f(\mathbf{b},\mathbf{u}) & f(\mathbf{b},\mathbf{v})
			\end{pmatrix}
			 \text{ and }
			\begin{pmatrix}
				t(\mathbf{a},\mathbf{u}) & t(\mathbf{a},\mathbf{v}) \\ t(\mathbf{b},\mathbf{u}) & t(\mathbf{b},\mathbf{v})
			\end{pmatrix}
		\end{equation*}
	\end{center}	
		where $f, t \in \Pol_{n+m}(\mathbf{A})$, $\mathbf{a}, \mathbf{b} \in A^n$, $\mathbf{u}, \mathbf{v} \in A^m$ with $\mathbf{a} \equiv_{\alpha} \mathbf{b}$ and $\mathbf{u} \equiv_{\beta} \mathbf{v}$, we have that:
	\begin{center}
		if $f(\mathbf{a},\mathbf{u}) \equiv_{\delta} t(\mathbf{a},\mathbf{u})$, $f(\mathbf{a},\mathbf{v}) \equiv_{\delta} t(\mathbf{a},\mathbf{v})$, and $f(\mathbf{b},\mathbf{u}) \equiv_{\delta} t(\mathbf{b},\mathbf{u})$ then $f(\mathbf{b},\mathbf{v}) \equiv_{\delta} t(\mathbf{b},\mathbf{v})$.
	\end{center}
	 We denote this ternary relation as $C_{2T}(\alpha,\beta;\delta)$ and refer to it as the \emph{two-term centralizer}. In our work, we utilize the TC-commutator introduced by Freese and McKenzie in \cite{FM.CTFC}. For properties and fundamental results concerning the commutator, we refer the reader to \cite{KK.TSOC}. Additionally, in Theorem \ref{ThmweakAlg}, we make use of the concepts of \emph{linear commutator} and \emph{two-term commutator}, denoted as $[\alpha,\beta]_L$ and $[\alpha,\beta]_{2T}$ respectively. These definitions can be found, for instance, in \cite{KK.TBTC} and \cite{Lip.ACOV}.
	
	Furthermore, a \emph{weak difference term} for a variety $\vv V$ is a ternary term $d(x, y, z)$ such that that for all  $\alg A \in \vv V$, $a, b \in A$, and $\theta \in \con(\mathbf{A})$ with $(a,b) \in \theta$, we have:
	
	\begin{equation*}
		d(b, b, a) \mathrel{[\theta, \theta]} a \mathrel{[\theta, \theta]}d(a, b, b),
	\end{equation*}

	For a comprehensive treatment of the notion, refer to \cite[Chapter $6$]{KK.TSOC}. We can observe that a weak difference term is Mal'cev on any block of an abelian congruence. Consequently, any abelian algebra within a variety with a weak difference term is affine \cite{KK.TSOC}. As a generalization of this definition we introduce the notion of \emph{commutator equation}. Let $p$ and $q$ terms for a language $\vv L$. Then we denote by $p \approx q$ an equation involving the two terms.  Throughout the paper, we refer to this type of equation as a \emph{term equation}. As a relaxation of this concept, we define the commutator equation. Let $\alg A$ be an algebra, and let $p$ and $q$ be $n$-ary terms of $\alg A$. We say that $\alg A$ \emph{satisfies} the commutator equation $p \approx_{C} q$, denoted as $\alg A \models p \approx_{C} q$, if for all $\theta \in \con(\alg A)$ and for all $a_1, \dots, a_n \in A$ in the same $\theta$-class, we have:
	%see \cite[Chapter $6$]{KK.TSOC} for a deep treatment of the notion. We can observe that a weak difference term is Mal’cev on any block of an abelian congruence, thus in particular any abelian algebra in a variety which has a weak difference term is affine, see \cite{KK.TSOC}. As a generalization of this definition we introduce the notion of \emph{commutator equation}. Let $p$ and $q$ terms for a language $\vv L$. Then we denote by $p \approx q$ an equation involving the two terms. Through all the paper we will refer to this standard type of equations as \emph{term equations}. As a relaxation of this concept we introduce the definition of commutator equation. Let $\alg A$ be an algebra and let $p$ and $q$ be $n$-ary terms of $\alg A$. We say that $\alg A$ \emph{satisfies} the commutator equation $p \approx_{C} q$, in symbols $A \models p \approx_{C} q$, if for all $\theta \in \con(\alg A)$ and for all $a_1, \dots a_n \in A$ in the same $\theta$-class, we have:
	\begin{equation*}
		p(a_1, \dots a_n) \mathrel{[\theta, \theta]} q(a_1, \dots a_n).
	\end{equation*}
	Clearly this concept of equation is weaker than the standard term equation. In fact, if an algebra satisfies the term equation $p \approx q$, then it also satisfies the commutator equation $p \approx_{C} q$
	
	Next, we recall the definition of \emph{Taylor term}, introduced in \cite{Tay.VOHL}, as it will be referenced extensively throughout this paper. An $n$-ary  \emph{Taylor term} is a term satisfying the following term equations: 
	\begin{align*}
		f(x_{11},\dots,x_{1n}) &\approx f(y_{11},\dots,y_{1n})
		\\&\ \cdot
		\\&\ \cdot
		\\&\ \cdot
		\\f(x_{n1},\dots,x_{nn}) &\approx f(y_{n1},\dots,y_{nn})
	\end{align*}
	with $x_{ij}, y_{ij} \in \{x,y\}$ and $x_{ii}\not=y_{ii}$, for all $i \in \{1, \dots, n\}$. These equations describe the class of varieties satisfying a non-trivial idempotent Mal'cev condition \cite{Tay.VOHL}, commonly referred to as  \emph{Taylor varieties}.
	
	Before moving forward, we recall the definition of \emph{Herringbone terms}, which can also be found in \cite[Section $8.2$]{KK.TSOC} or in \cite{Lip.ACOV}.  These terms are of particular interest due to their significant connection with well-behaved commutators in varieties, as described in \cite{KK.TBTC} and \cite{Lip.ACOV}.
	\begin{Def}\label{DefHerring}
		We call $\emph{Herringbone terms}$ the family of lattice terms $\{y^i\}_{i \in \NN}$% and $\{z^i\}_{i \in \NN}$  
    in the variables $\{x,y,z\}$ defined by:
		\begin{align*}
			&y^{0}(x,y,z) = y; 
   %z^{0}(x,y,z) = z;
			\\&y^{n+1}(x,y,z) = y \vee (x \wedge y^{n}(x,z,y)).	
			%\\z^{n+1}(x,y,z) &= z \vee (x \wedge y^{n}(x,y,z)).
		\end{align*}
	\end{Def}

	We adopt the superscript notation instead of the conventional underscript notation, as the latter conflicts with the customary use of  underscript for denoting sequences of variables. It is worth noting that these lattice terms give rise to a potentially infinite ascending chain and that are often presented in their dual version. Additionally, let $p$ be a term for the language ${\wedge,\vee,\circ}$, and consider $k \in \mathbb{N}$. We denote by $p^{(k)}$ the $\{\wedge,\circ\}$-term obtained by replacing every occurrence of $\vee$ with the $k$-fold relational product $\circ^{(k)}$.
	%We use the superscript instead of the standard subscript since the latter is in conflict to the usual subscrit notation for a sequence of variables. Note that this lattice terms form two possibly infinite ascending chains. Furthermore, let $p$ be a term for the language $\{\wedge,\vee,\circ\}$. Let $k \in \NN$. We denote by $p^{(k)}$ the $\{\wedge,\circ\}$-term obtained from $p$ substituting any occurrence of $\vee$ with the $k$-fold relational product $\circ^{(k)}$. %We denote by \index{$[n]$}$[n]$ the set $\{i \in \NN\mid 1 \leq i \leq n\}$.
	\section{Labelled graphs and regular terms}
	\label{LabelledG}
	To present the primary findings of Section \ref{Alg}, we recall the definition \emph{labelled graph associated with a term} as usually defined 
 (cf. e.g. \cite{Cze.AMTC, CD.HSWW, KK.TSOC}). Let us begin by providing the definition of a labelled graph.
	
	\begin{Def}
		Let $S$ be a set of labels. Then a \emph{labelled graph} is a directed graph $\mathbf{G} = (V_{\mathbf{G}},E_{\mathbf{G}})$ with a labelling function $l: E_{\mathbf{G}} \rightarrow S$.
	\end{Def}

	Following \cite{KK.TSOC}, let $p$ be a $\{\wedge,\circ\}$-term in the variables $\{X_1, \dots, X_t\}$. We present a construction that yields a finite sequence of labelled graphs  $\{\mathbf{G}_i(p)\}_{i \in I}$. This sequence starts with the labelled graph $\mathbf{G}_1(p) = (\{x_1,x_2\},$ $ \{(x_1,x_2)\})$ with $l((x_1,x_2)) = p$ having an edge $(x_1,x_2)$ labelled with $p$, connecting two vertices $x_1$ and $x_2$. For $s \geq 1$, we proceed from the labelled graph $\mathbf{G}_s(p)$ with vertices $V_{\mathbf{G}_s(p)} = \{x_1, \dots, x_{k}\}$ by selecting (if possible) $w  \not\in \{X_1, \dots, X_t\}$ such that $w$ is a label of an edge $(x_i, x_j)$ in $\mathbf{G}_s(p)$. Then we have two cases.
	 %Then from the labelled graph $\mathbf{G}_s(p)$ we continue selecting $w \not\in \{x_1, \dots, x_t\}$ such that $w$ is a label of an edge $(y_i, y_j)$ of $\mathbf{G}_s(p)$.  Then we have two cases.
	\begin{itemize}
		\item If $w = u \wedge v$ for some $u$ and $v$ $\{\wedge,\circ\}$-terms. Then $\mathbf{G}_{s+1}(p)$ is obtained from $\mathbf{G}_s(p)$ by replacing the edge $(x_i,x_j)$ labelled $w$ with two edges connecting the same vertices $x_i$ and $x_j$ labelled $u$ and $v$ respectively;
		\item if $w = u \circ v$ for some $u$ and $v$ $\{\wedge,\circ\}$-terms. Then $\mathbf{G}_{s+1}(p)$  is obtained from $\mathbf{G}_s(p)$ by introducing a new vertex $x_{k+1}$ and replacing the edge $(x_i,x_j)$ labelled $w$ with two edges $(x_i,x_{k+1})$ and $(x_{k+1},x_j)$, labelled $u$ and $v$ respectively, connecting the vertices in serial through $x_{k+1}$.
	\end{itemize} 
	The construction ends when for some $n \in \NN$  none of the above steps can be executed for the labelled graph $\mathbf{G}_{n}(p)$. Thus, when $l(e) \in  \{X_1, \dots, X_t\}$ for every $e$ edge in $\mathbf{G}_{n}(p)$.
	
	It is worth noting that the choice of $w$ in every step may lead to different sequences derived from a term $p$. However, we observe that the last graph in the sequence remains the same up to reordering of the vertices. We call the last graph of the sequence $\{\mathbf{G}_i(p)\}_{i \in I}$  \emph{labelled graph associated with} $p$, denoted by $\mathbf{G}(p)$. 
	The main reason to introduce $\mathbf{G}(p)$ is stated ind \cite[Proposition 3.1]{CD.HSWW} and in Claim $4.8$ of \cite{KK.TSOC}. The latter can be generalized to tolerances as shown in \cite[Proposition $2.1$]{Lip.FCIT} and even to relations in general \cite[Proposition $3.3$]{Fio.MCCT}. We include the most general version of this technical result in order to use it in Section \ref{Alg}.
	
	\begin{Prop}[Proposition $3.3$ of \cite{Fio.MCCT}]
		\label{PropKK}
		Let $\mathbf{A}$ be an algebra, let $p$ be a $\{\circ ,\wedge\}$-term in the variables $\{X_m\}_{m \in I}$, and let $R_i \subseteq A \times A$ for $i \in I$. Then:
		
		\begin{enumerate}
			\item[(1)]  Let $V_{\mathbf{G}(p)} \rightarrow A$: $x_s \mapsto a_s$ be an assignment such that for all edges $(x_i, x_j)$  with label $X_k$ of $\mathbf{G}(p)$ and $k \in I$, we have $(a_i,a_j) \in R_k$. Then $(a_1, a_2) \in p(R_1,\dots,R_n)$.
			
			\item[(2)] Conversely, given any $(a_1, a_2) \in p(R_1,\dots,R_n)$, there is an assignment $V_{\mathbf{G}(p)} \rightarrow A$: $x_s \mapsto a_s$ extending $x_1 \mapsto a_1$, $x_2 \mapsto a_2$ such that $(a_i,a_j) \in R_k$ whenever $(x_i, x_j)$ is an edge labelled with $X_k$ of $\mathbf{G}(p)$ and $k \in I$, where $(x_1,x_2)$ is the only edge of the graph $\mathbf{G}_1(p)$.
		\end{enumerate}
	\end{Prop}

	\section{Taylor varieties and abelian algebras}
	\label{Taylor}
	
	In this section our aim is to show a property describing Taylor varieties logically weaker then the standard presented in \cite{Tay.VOHL}. Taylor varieties have attracted significant attention for several reasons. Surprisingly, the class of Taylor varieties constitutes a strong Mal'cev class \cite{Ols.TWNI}, and this particular class of varieties exhibits a profound connection with the Feder-Vardi conjecture, which was independently proven to be true in \cite{Bul.ADTF} and \cite{Zhuk.APOT}.  Furthermore, contemporary research has uncovered insightful findings regarding the relationship between commutators and Taylor varieties \cite{Ker.RMDO}. Thus, Taylor varieties have consistently represented a central theme of investigation in the field of universal algebra.

    Before delving deeply into the theory of commutator equations, we begin with a straightforward observation.

	\begin{Rem}\label{idempRem}
		Let $\vv V$ be a variety satisfying the commutator equation $x \approx_{C} p(x,\dots, x)$. Then $p$ is idempotent.
	\end{Rem}
	
	The remark follows from the fact that the commutator equation must hold also for the $0$ congruence.
	
	In order to prove the main result of the section we present a lemma that can be seen in \cite{Lip.ACOV} in a slightly different version without a proof. From now on when not specified $\beta^n = y^n(\alpha, \beta, \gamma)$ and $\gamma^n = y^n(\alpha, \gamma, \beta)$.
	
	\begin{Lem}\label{herringLem}
		Let $\alg{A}$ be an algebra and $\alpha, \beta, \gamma \in \con(\alg A)$. Let $\delta = \bigcup_{n \in \mathbb{N}}(\alpha \wedge \beta^{n})$. Then,
		\begin{equation*}
			[\alpha \wedge (\beta \vee \gamma), \alpha \wedge (\beta \vee \gamma)] \leq \delta.
		\end{equation*}
	\end{Lem}
	
	\begin{proof}
		We observe that Definition \ref{DefHerring} yields that the chains of congruences $\{\alpha \wedge \beta^{n}\}_{n \in \mathbb{N}}$ and $\{\alpha \wedge \gamma^{n}\}_{n \in \mathbb{N}}$ are cofinal since:
		\begin{equation*}
			\alpha \wedge \beta^{n+1} \geq \alpha \wedge \gamma^{n} \geq \alpha \wedge \beta^{n-1}.
		\end{equation*}
		Hence $$\delta = \bigcup_{n \in \mathbb{N}}(\alpha \wedge \beta^{n})=\bigcup_{n \in \mathbb{N}}(\alpha \wedge \gamma^{n}).$$ This identity, the equation $\beta^{n+1} = \beta \vee (\a \wedge \g^n)$ , and the fact that $\wedge$ distributes over the union implies:
		\begin{equation}\label{cofinaleq}
			\alpha \wedge (\beta \vee (\alpha \wedge \delta))= \alpha \wedge (\beta \vee (\alpha \wedge \bigcup_{n \in \mathbb{N}}(\alpha \wedge \gamma^{n})) = \bigcup_{n \in \mathbb{N}}(\alpha \wedge \beta^{n}) = \delta.
		\end{equation}
		Clearly a similar equality holds substituting $\beta$ with $\gamma$. From this two identities and the properties centralizer relation, see \cite[Theorem $2.19$ $(8)$]{KK.TSOC}, we obtain that $C(\beta,\alpha;\delta)$ and $C(\gamma,\alpha;\delta)$. The semidistributivity of the centralizer in the first component (\cite[Theorem $2.19$ $(5)$]{KK.TSOC}) yields $C(\beta \vee \gamma,\alpha;\delta)$. From the definition of commutator we obtain $[\beta \vee \gamma,\alpha] \leq \delta$ and the claim follows from the monotonicity of the commutator.
	\end{proof}
	
	Lemma \ref{herringLem}, in conjunction with the following theorem, plays a crucial role in proving the main result of this section. The theorem establishes a connection between Ol\v{s}'{a}k's result \cite{Ols.TWNI} and the equation that characterizes Taylor varieties in \cite{KK.TBTC}. By leveraging \cite{Ols.TWNI} and the congruence equation in \cite[Lemma $4.6$]{KK.TBTC}, we obtain a simplified version of this equation. The original congruence equation characterizing Taylor varieties in \cite{KK.TBTC} is derived using an arbitrary $n$-ary Taylor term, but with the introduction of a $6$-ary Ol\v{s}'{a}k term, we can refine it to a more comprehensible form that involves a fixed number of congruences.
	
	\begin{Thm}\label{TaylorThmorig}
		Let $\vv{V}$ be a variety. Then the following are equivalent:	
		\begin{enumerate}		
			\item[(1)] $\vv{V}$ is Taylor;
			\item[(2)] $\vv{V}$ satisfies the term equations:
			\begin{equation}\label{Olsak eq}
				t(x,y,y,y,x,x) \approx  t(y,x,y,x,y,x) \approx 	t(y,y, x,x,x,y);
			\end{equation}
			for some idempotent term $t$.
			\item[(3)] $\vv{V}$ satisfies the following congruence equation in the variables $\{\alpha_1, \dots, \alpha_6, \beta_1, \dots, \beta_6\}$:
			\begin{equation}\label{Tayeq}
				\bigwedge^{6}_{i = 1}(\a_{i} \circ \b_{i}) \leq (\bigvee^{6}_{i=1}\a_{i} \wedge \bigwedge^{2}_{i=1} (\g \vee \theta_{i})) \vee (\bigvee^{6}_{i=1}\b_{i} \wedge \bigwedge^{2}_{i=1} (\g \vee \theta_{i}));
			\end{equation}
			where $\gamma =  \bigwedge^{6}_{i = 1} (\a_{i} \vee \beta_i)$ and 
			\begin{equation*}
				\theta_{i}= (\bigvee_{j \in L_{i}} \a_{j} \vee \bigvee_{j \in L'_{i}} \b_{j}) \wedge (\bigvee_{j \in R_{i}} \a_{j} \vee \bigvee_{j \in R'_{i}} \b_{j})
			\end{equation*}
			with 
			\begin{align*}
				&L_1 = L_2 =\{1,5,6\}  &L'_{1} = L'_2 = \{ 2,3,4 \}
				\\&R_{1} = \{ 2,4,6\} &R'_{1} = \{1,3,5 \}
				\\&R_{2} = \{ 3,4,5\} &R'_{2} = \{1,2,6 \}
			\end{align*}
		\end{enumerate}
	\end{Thm}
	
	\begin{proof}
		The equivalence of $(1)$ and $(2)$ follows from Ol\v{s}\'{a}k's result \cite{Ols.TWNI}. $(2) \Rightarrow (3)$ follows from \cite[Lemma $4.6$]{KK.TBTC} adjusting the proof for the $6$-ary Taylor term in $(2)$. For sake of completeness we include the modified version of the proof. Let $\alg A \in \vv V$ and let $(a,b) \in \bigwedge^{6}_{i = 1}(\a_{i} \circ \b_{i})$. Thus there exist $u_1,\dots u_6 \in A$ such that $a \mathrel{\a_i} u_i  \mathrel{\b_i} b$, for all $i \in \{1,\dots,6\}$. Let $u = t^{\alg A}(u_1,\dots,u_6)$ and $v = t^{\alg A}(a,b,b,b,a,a) = t^{\alg A}(b,a,b,a,b,a) = t^{\alg A}(b,b,a,a,a,b)$. We prove that:
		\begin{enumerate}
			\item[(a)] $(a,u) \in \bigvee^{6}_{i=1}\a_{i}$;
			\item[(b)] $(u,b) \in \bigvee^{6}_{i=1}\b_{i}$;
			\item[(c)] $ (a,u), (u,b) \in \bigwedge^{2}_{i=1} (\g \circ \theta_{i})$.
		\end{enumerate}
		For claims (a) and (b) we observe that:
		\begin{align*}
			a &=  t^{\alg A}(a,\dots, a) \mathrel{\a_1} \circ \cdots \circ \mathrel{\a_6}  t^{\alg A}(u_1,\dots,u_6)
			\\b &= t^{\alg A}(b,\dots, b) \mathrel{\b_1} \circ \cdots  \circ \mathrel{\b_6}  t^{\alg A}(u_1,\dots,u_6)
		\end{align*}
		For claim (c) we notice that $a$, $b$ and $v$ are in the same $\g$-class. Furthermore:
		\begin{equation*}
			u = t^{\alg A}(u_1,\dots,u_n)  \mathrel{\theta_{i}} v.
		\end{equation*}
		Thus $a \mathrel{\gamma} v \mathrel{\theta_{i}} u$ and $b \mathrel{\gamma} v \mathrel{\theta_{i}} u$ and (c) holds. Putting (a), (b), and (c) together we obtain that $(a,b)$ is in the right side of the inequality (\ref{Tayeq}) and $(3)$ holds.
		For $(3) \Rightarrow (1)$ we can observe that the equation (\ref{Tayeq}) is a spacial case of the equation in \cite[proof of  Lemma $4.6$]{KK.TBTC}, hence (\ref{Tayeq}) is non-trivial and thus it implies the satisfaction of a non-trivial idempotent Mal'cev condition. This can be seen considering an $8$-element set $\{a, b, u_1, \dots, u_n\}$ with the partitions $\a_i$ generated by $(a,u_i)$ and  $\b_i$ generated by $(b,u_i)$, for all $i \in \{1,\dots,6\}$. Hence $(1)$, $(2)$ and $(3)$ are equivalent.
	\end{proof}

	It is worth noting that equation (\ref{Tayeq}) can be further simplified by utilizing the relation product and entirely avoiding the use of $\join$. This approach leads to the generation of a strong Mal'cev condition through the Pixley-Wille algorithm \cite{Pix.LMC,Wil.K}. With these considerations in place, we are ready to prove the main result of the section which uses Theorem \ref{TaylorThmorig} to find a Mal'cev condition equivalent to be Taylor involving commutator equations. 
	
	\begin{Thm}\label{TaylorThm}
		Let $\vv{V}$ be a variety. Then the following are equivalent:	
		\begin{enumerate}		
			\item[(1)] $\vv{V}$ is Taylor;
			\item[(2)] the variety $\vv{V}'$ generated by abelian algebras of the idempotent reduct of $\vv{V}$ is Taylor;
			
			\item[(3)] $\vv V$ satisfies the commutator equations:
			\begin{align*}\label{weakTay}
				t(x,y,y,y,x,x) &\approx_{C}  t(y,x,y,x,y,x) \approx_{C} t(y,y, x,x,x,y)
				\\t(x,x,x,x,x,x) &\approx_{C} x;
			\end{align*}
			
			\item[(4)] $\vv{V}$ satisfies the following congruence equation:
			\begin{equation*}\label{Tayeqweak}
				\bigwedge^{6}_{i = 1}(\a_{i} \circ \b_{i}) \leq (\bigvee^{6}_{i=1}\a_{i} \wedge \bigwedge^{2}_{i=1} (\tau \vee \theta_{i})) \vee (\bigvee^{6}_{i=1}\b_{i} \wedge \bigwedge^{2}_{i=1} (\tau \vee \theta_{i}));
			\end{equation*}
			where 
			\begin{equation*}
				\theta_{i}=(\bigvee_{j \in L_{i}} \a_{j} \vee \bigvee_{j \in L'_{i}} \b_{j}) \wedge ( [\tau,\tau] \circ \bigvee_{j \in R_{i}} \a_{j} \vee \bigvee_{j \in R'_{i}} \b_{j})
			\end{equation*}
			with $\tau = \bigwedge^{6}_{i = 1}(\a_{i} \vee\b_{i})$ and
			\begin{align*}
				&L_1 = L_2 =\{1,5,6\}  &L'_{1} = L'_2 = \{ 2,3,4 \}
				\\&R_{1} = \{ 2,4,6\} &R'_{1} = \{1,3,5 \}
				\\&R_{2} = \{ 3,4,5\} &R'_{2} = \{1,2,6 \}
			\end{align*}
			\item[(5)]  there exists $n \in \NN$ such that $\vv V$ satisfies the following congruence equation in the variables $\{\alpha_1, \dots, \alpha_6, \beta_1, \dots, \beta_6\}$:
			\begin{equation}\label{weakTayHarreq}
				\bigwedge^{6}_{i = 1}(\a_{i} \circ \b_{i}) \leq (\bigvee^{6}_{i=1}\a_{i} \wedge \bigwedge^{2}_{i=1} (\tau \vee \theta_{i})) \vee (\bigvee^{6}_{i=1}\b_{i} \wedge \bigwedge^{2}_{i=1} (\tau \vee \theta_{i}));
			\end{equation}
			where $\tau =  \bigwedge^{6}_{i = 1} (\a_{i} \vee \beta_i)$ and 
			\begin{equation*}
				\theta_{i}=(\bigvee_{j \in L_{i}} \a_{j} \vee \bigvee_{j \in L'_{i}} \b_{j}) \wedge ((\bigwedge_{i=1}^{5}(\a_i \vee \b_i) \wedge \beta_6^n) \circ \bigvee_{j \in R_{i}} \a_{j} \vee \bigvee_{j \in R'_{i}} \b_{j})
			\end{equation*}
			with $\beta_6^n = y^n(\bigwedge^{5}_{i = 1} (\a_{i} \vee \beta_i), \beta_6, \alpha_6)$ and
			\begin{align*}
				&L_1 = L_2 =\{1,5,6\}  &L'_{1} = L'_2 = \{ 2,3,4 \}
				\\&R_{1} = \{ 2,4,6\} &R'_{1} = \{1,3,5 \}
				\\&R_{2} = \{ 3,4,5\} &R'_{2} = \{1,2,6 \}
			\end{align*}
		\end{enumerate}
	\end{Thm}
	
	\begin{proof}
		$(1) \Rightarrow (2)$ is trivial since if $\vv V$ is Taylor then also the idempotent reduct of $\vv V$ is Taylor.
		
		For $(2) \Rightarrow (3)$, let us assume that the variety $\vv{V}'$ generated by abelian algebras of the idempotent reduct of $\vv{V}$ is Taylor. Thus it has an  Ol\v{s}\'{a}k term $t$ that satisfies the equations (\ref{Olsak eq}) in Theorem \ref{TaylorThmorig}.  Let $\alg A \in \vv V$ and let $\theta \in \con(\alg A)$.  Factoring $\alg A$ by $[\theta, \theta]$ we see that is sufficient to verify the equations in (3) for abelian congruences. Let $\alg A'$ be the idempotent reduct of $\alg A$. We can observe that every $\theta$-class of $\alg A$ is an abelian subalgebra of  $\alg A'$ and thus it has a term $t$ satisfying (\ref{Olsak eq}) by (2). We claim that $t$ satisfies the equations in (3) for $\vv V$. Let $(a,b) \in \theta$, then $t^{\alg A'}(a,b,b,b,a,a) = t^{\alg A'}(b,a,b,a,b,a) = t^{\alg A'}(b,b,a,a,a,b)$ and $t^{\alg A'}(a,a,a,a,a,a) = a$.
		The commutator equations (\ref{Olsak eq}) follow directly remembering that a quotient by $[\theta, \theta]$ has been applied.
		
		For $(3) \Rightarrow (4)$, suppose that $(3)$ holds. Then let $\alg A \in \vv V$ and let $\a_1, \dots, a_6, \b_1, \dots, \b_6 \in \con (\alg A)$ with $(a,b) \in \bigwedge^{6}_{i = 1}(\a_{i} \circ \b_{i})$. Then, with a similar argument of $(2) \Rightarrow (3)$ in the proof of  Theorem \ref{TaylorThmorig} we have that there exist $u_1,\dots u_6 \in A$ such that $a \mathrel{\a_i} u_i  \mathrel{\b_i} b$, for all $i \in \{1,\dots,6\}$. Let $u = t^{\alg A}(u_1,\dots,u_6)$ and $v = t^{\alg A}(a,b,b,b,a,a)$. Since $(a,b) \in \tau$ we have that $v \mathrel{[\tau, \tau]} t^{\alg A}(b,a,b,a,b,a)  \mathrel{[\tau, \tau]} t^{\alg A}(b,b,a,a,a,b)$, by the equations in (3). We prove that:
		\begin{enumerate}
			\item[(a)] $(a,u) \in \bigvee^{6}_{i=1}\a_{i}$;
			\item[(b)] $(u,b) \in \bigvee^{6}_{i=1}\b_{i}$;
			\item[(c)] $ (a,u), (u,b) \in \bigwedge^{2}_{i=1} (\tau \circ \theta_{i})$.
		\end{enumerate}
		For claims (a) and (b) we can observe that $t$ is idempotent, by Remark \ref{idempRem}, and thus
		\begin{align*}
			a &=  t^{\alg A}(a,\dots, a) \mathrel{\a_1} \circ \cdots \circ \mathrel{\a_6}  t^{\alg A}(u_1,\dots,u_6)
			\\b &= t^{\alg A}(b,\dots, b) \mathrel{\b_1} \circ \cdots  \circ \mathrel{\b_6}  t^{\alg A}(u_1,\dots,u_6).
		\end{align*}
		For claim (c) we notice that $a$, $b$ and $v$ are in the same $\tau$-class. Furthermore, we can prove that $v \mathrel{\theta_{i}} u$, for $i = 1,2$.
		\begin{align*}
			t^{\alg A}(u_1,\dots,u_n)  & \mathrel{\bigvee_{j \in L_{1}} \a_{j} \vee \bigvee_{j \in L'_{1}} \b_{j}}  v
			\\t^{\alg A}(u_1,\dots,u_n) & \mathrel{\bigvee_{j \in R_{1}} \a_{j} \vee \bigvee_{j \in R'_{1}} \b_{j}} t^{\alg A}(b,a,b,a,b,a) \mathrel{[\tau,\tau]} v
			\\t^{\alg A}(u_1,\dots,u_n) & \mathrel{\bigvee_{j \in R_{2}} \a_{j} \vee \bigvee_{j \in R'_{2}} \b_{j}} t^{\alg A}(b,b,a,a,a,b) \mathrel{[\tau,\tau]}  v.
		\end{align*}
		Hence,  $a \mathrel{\tau} v \mathrel{\theta_{i}} u$ and $b \mathrel{\tau} v \mathrel{\theta_{i}} u$ thus (c) holds. Putting (a), (b), and (c) together we obtain that $(a,b)$ is in the right side of the inequality in $(4)$.
		
		$(4) \Rightarrow (5)$ follows applying Lemma \ref{herringLem} to the equation in $(4)$ with $\alpha = \bigwedge^{5}_{i=1}(\alpha_i \vee \beta_i)$, $\beta = \beta_6$, and $\gamma = \alpha_6$. Thus $[\bigwedge^{6}_{i=1}(\alpha_i \vee \beta_i), \bigwedge^{6}_{i=1}(\alpha_i \vee \beta_i)] \leq \bigcup_{n \in \mathbb{N}}(\alpha \wedge \b^n)$ and hence there exists $n \in \NN$ such that:
		\begin{equation*}
				[\tau,\tau] \leq \bigwedge_{i=1}^{5}(\a_i \vee \b_i) \wedge \beta_6^n.
		\end{equation*}
		This yields the equation in $(5)$.
		For $(5) \Rightarrow (1)$ we show that equation (\ref{Tayeqweak}) is non-trivial and this implies $(1)$ through the Pixley-Wille algorithm, see \cite{Pix.LMC, Wil.K}.  Let $\alg A = \langle \{a,b, c_1 \dots c_6\}, \pi \rangle$ be an algebra whose basic operations are only projections. Clearly every equivalence class is a congruence in $\con(\alg A)$. Let $\alpha_i$ be the partitions which identify only $\{a, c_i\}$ and let $\b_i$ be the partitions which identify only $\{b, c_i\}$, for all $i \in \{1,\dots,6\}$.  Then $(a,b)$ is in the left side of (\ref{weakTayHarreq}). We prove that $(a,b)$ is not in the right side. First we can observe that  $\bigwedge_{i=1}^{5}(\a_i \vee \b_i) \wedge \b_6^n = 0$ for all $n \in \NN$, since $\b_6^n = \b_6$ for all $n \in \NN$. Furthermore, $\bigvee_{i=1}^{6}\a_i$ has exactly the two classes $\{a, c_1,\dots, c_6\}$ and $\{b\}$ while the two classes of  $\bigvee_{i=1}^{6}\b_i$ are $\{c_1,\dots, c_6, b\}$ and $\{a\}$. Moreover, we can observe that $(a,c_j), (c_j,b) \notin \bigwedge^{2}_{i=1} \tau \vee \theta_i$ for all $j = 1,\dots, 6$. Consequently $\bigwedge_{i=1}^2 \tau \vee \theta_i $ has $\{a,b\}$ as the only non-trivial class. Hence we can conclude that $\{a\}$ is a class of $\bigvee^{6}_{i=1}\a_{i} \wedge \bigwedge^{2}_{i=1} (\tau \vee \theta_{i})$   and  $\{b\}$ is a class of $\bigvee^{6}_{i=1}\b_{i} \wedge \bigwedge^{2}_{i=1} (\tau \vee \theta_{i})$. Thus $(a,b) \notin (\bigvee^{6}_{i=1}\a_{i} \wedge \bigwedge^{2}_{i=1} (\tau \vee \theta_{i})) \vee (\bigvee^{6}_{i=1}\b_{i} \wedge \bigwedge^{2}_{i=1} (\tau \vee \theta_{i}))$ as wanted.
	\end{proof}
	
	The previous theorem is an improvement of  Taylor's  characterization of the class of varieties satisfying a non-trivial idempotent Mal'cev condition obtained with a relaxation of the assumptions. It is noteworthy that the proof does not rely exclusively on the Ol\v{s}'{a}k term and can be applied using a generic $n$-ary Taylor term instead. However, employing the Ol\v{s}'{a}k term enables the derivation of simpler equations. The previous theorem can also be seen as a consequence of \cite[Theorem 3.13]{KK.TSOC}, which shows that Taylor varieties omit strongly abelian congruences, a stronger result. Nonetheless, we have included the previous result in our discussion due to the applicability of the technique employed in its proof within a broader context. In the subsequent section, we delve further this intuition and explore its implications.
	
	\section{A Pixley-Wille type algorithm for commutator equations}
	\label{Alg}

	In this section our aim is to develop a Pixley-Wille type algorithm for commutator equations, something that can be partially seen in \cite[Theorem $3.4$]{Lip.ACOV} without a proof. Investigating this problem is of interest since a systematic way of characterizing  Mal'cev conditions described by commutator equations can produce idempotent Mal'cev conditions logically weaker than their congruence counterpart.
	
	Let $p \leq q$ be an equation for the language $\{\wedge, \circ\}$ in the variables $\{X_s\}_{s \in I}$. Let $\mathbf{G}(p)$ and $\mathbf{G}(q)$ be the graph associated to respectively $p$ and $q$ obtained through the procedure described in Section \ref{LabelledG} and let $\{x_1, \dots, x_n\}$ be the set of vertices of $\mathbf{G}(p)$. For all $s \in I$, we define:
	\begin{align}
		\label{eqp}
		T_s(p) &:= \{(x_i,x_j) \mid (x_i,x_j) \text{ is an edge of } \mathbf{G}(p) \text{ with label } X_s \}
		\\Tt^n(q) &:=\{(t_i,t_j) \mid (x_i,x_j) \text{ is an edge of } \mathbf{G}(q)\}\nonumber
		\\Tt_s^n(q) &:=\{(t_i,t_j) \mid (x_i,x_j) \text{ is an edge of } \mathbf{G}(q) \text{ with label } X_s\},\nonumber
	\end{align}
	where the elements $\{t_1,\dots,t_l\}$ of the pairs in $Tt^n(q)$ are $n$-ary terms with $t_1 = x_1$ and $t_2 = x_n$ in the variables $\{x_1,\dots,x_n\}$. Let $R \subseteq A \times A$ be a relation over the set $A$. We denote by $\Eqv(R)$ the equivalence relation generated by $R$. We define $\Eq_s(p \leq q)$ as the set of all equations of the form:
	\begin{equation}\label{Pixeq}
		t_i(\sigma_s(x_{1}),\dots,\sigma_s(x_{m})) \approx t_j(\sigma_s(x_{1}),\dots,\sigma_s(x_{m}))
	\end{equation}
	such that $(t_i,t_j) \in Tt^n_s(q)$, where $n$ is the number of vertices in $\mathbf{G}(p)$ and $\sigma_s$ is a transversal of $\Eqv(T_s(p))$, i.e. a map collapsing variables in the same $\Eqv(T_s(p))$-class into one fixed element of the class. %only the indices of variables in a same class of $\Eqv(T_s(p))$. %This means that the variables that are in the equivalence relation generated by the pairs in $T_s(p)$ are collapsed. 
    Furthermore, we define $\Eq(p \leq q) = \bigcup_{s \in I}\Eq_s(p \leq q)$ which are exactly the equations produced by the Pixley-Wille algorithm \cite{Pix.LMC,Wil.K} characterizing the Mal'cev condition described by the congruence equation $p \leq q$.
	
	%In the previous definition we  introduced terms of unspecified type. We consider those terms to be either of one letter $x_1, x_2, \dots$ or composed by a primitive operation symbols applied to letters, e.g. $t(x_1,\dots,$ $ x_n)$. Namely, those terms are placeholders used to write equations that can be instantiated in the language of a given variety. 
    
    %Let  $\Sigma$ be a set  terms equations of unspecified type and let $\vv{V}$ be a variety of type $\tau$. Then a $\tau$-\emph{realization} of $\Sigma$ is a set of equations obtained from $\Sigma$ by replacing each operation symbol, in all of its occurrences in $\Sigma$, by some fixed term symbol of type $\tau$. Furthermore, a variety $\vv{V}$ satisfies a set of equations $\Sigma$ of terms of unspecified type if there exists a $\tau$-realization of $\Sigma$ such that the resulting equations hold in $\vv{V}$. Henceforth, in accordance with common practice in the literature, we will use the same symbols to represent both the terms used as placeholders in equations and their $\tau$-realizations in a variety of type $\tau$.
	
	Let $p$ be a term for the language $\{\wedge, \circ\}$ and $q$ be a term for the language $\{\wedge, \circ, \vee\}$. We say that a variety $\vv{V}$ satisfies $\Eq(p \leq q)$ if there exists $k \in \NN$ such that $\vv{V}$ satisfies $\Eq(p \leq q^{(k)})$.
	
	We adapt the final part of the algorithm presented in \cite{Pix.LMC, Wil.K} to obtain a set of equations $\Eq_{C}(p \leq q)$ that characterizes the Mal'cev condition describing the class of varieties which satisfy $p \leq q$ over the algebras of the variety generated by abelian algebras of the idempotent reduct.
	\begin{Alg}
		\label{PixWilAlg}
		\theoremstyle{definition}
		Let $p \leq q$ be an $m$-ary equation for the language $\{\wedge, \circ\}$. Let $\mathbf{G}(p)$ and $\mathbf{G}(q)$ be obtained from $p$ and $q$ with the procedure in Section \ref{LabelledG}. For all $1 \leq s \leq m$, let us consider $T_s(p)$ and $Tt^n_s(q)$ as in \eqref{eqp}, where $n$ is the number of vertices of $\alg G(p)$. We define $\Eq^s_{C}(p \leq q)$ as the of the set of all equations of the form:
\begin{equation*}t_i(\sigma_s(x_{1}),\dots,\sigma_s(x_{n})) \approx_{C} t_j(\sigma_s(x_{1}),\dots,\sigma_s(x_{n}))
		\end{equation*}
		such that $(t_i,t_j) \in Tt^n_s(q)$ and $\sigma_s$ is a transversal of $\Eqv(T_s(p))$.
  Furthermore, we define $\Eq_{C}(p \leq q) = \bigcup_{1 \leq s \leq m} \Eq^s_{C}(p \leq q)$.
	\end{Alg}

	In order to prove the main result of the section we need the following technical Lemma inspired by \cite[Theorem $3.1$]{Lip.ACOV}.
	
	\begin{Lem}\label{lem:classher}
		Let $\vv{V}$ be a variety, let $\alg{F}_{\vv{V}}(3)$ be the $3$-generated free algebra of $\vv{V}$, let $\a = \Cg(\{(x,z)\}), \b = \Cg(\{(x,y)\}), \g  = \Cg(\{(y,z)\})$, and let $s,t \in \alg{F}_{\vv{V}}(3)$ with $(s, t) \in \a \wedge y^n(\alpha,\b,\g)$ for some $n \in \NN$. Then, for all $\alg{A} \in \vv{V}$ and $\delta \in \con(\alg{A})$ with $[\delta,\delta]_{2T} = 0$, we have $t^{\alg{A}}(a,b,c) = s^{\alg{A}}(a,b,c)$ for all $a,b,c$ in the same $\delta$-class.
	\end{Lem}

	\begin{proof}
		First we prove that for all $\alg{A} \in \vv{V}$, $\delta \in \con(\alg{A})$ with $[\delta,\delta]_{2T} = 0$, and $u, v$ ternary terms such that $u^{\alg{A}}(a,b,a) = v^{\alg{A}}(a,b,a)$ and $u^{\alg{A}}(a,b,b)$ $ = v^{\alg{A}}(a,b,b)$ for all $a,b$ in the same $\delta$-class, then $u^{\alg A} = v^{\alg A}$ over the same $\delta$-block. This claim is included in a symmetric version in \cite[Theorem $3.1$]{Lip.ACOV}. In order to prove the claim we can observe that for all $a,b,c$ in the same $\delta$-class, we can consider the matrices:
		\begin{center}
			\begin{equation*}
				\begin{pmatrix}
					u^{\alg{A}}(b,c,b) & 	u^{\alg{A}}(a,c,c) \\ 	u^{\alg{A}}(b,b,b) & 	u^{\alg{A}}(a,b,c)
				\end{pmatrix}
				\text{ and } 
				\begin{pmatrix}
					v^{\alg{A}}(b,c,b) & 	v^{\alg{A}}(a,c,c) \\ 	v^{\alg{A}}(b,b,b) & 	v^{\alg{A}}(a,b,c)
				\end{pmatrix}
				%\in M(\delta,\delta)
			\end{equation*}
		\end{center}
        By definition of two-term commutator we have
		 $u^{\alg A} \mathrel{[\delta, \delta]_{2T}} v^{\alg A}$ over the same $\delta$-block and this yields $u^{\alg A} = v^{\alg A}$ over the same $\delta$-class. Moreover, suppose that $(s, t) \in \a \wedge y^1(\alpha,\b,\g) = \a \wedge (\b \vee (\a \wedge \g))$. Then there exists $k \in \NN$ such that $(s, t) \in \a \wedge (\b \circ^{(2k+1)}  (\a  \wedge \g))$. Thus $s \mathrel{\a} t$ and there exist $v_1, w_1, \dots, v_n, w_n \in \alg{F}_{\vv{V}}(3)$ such that 
	\begin{equation*}
		s \mathrel{\b} v_1 \mathrel{\a \wedge \g} w_1 \cdots v_k \mathrel{\a \wedge \g} w_k \mathrel{\b} t.
	\end{equation*}
	By a standard Mal'cev argument (\cite[Lemma $12.1$]{BS.ACIU}), $s(x,y,x) = t(x,y,x)$,  $s(x,x,y) = v_{1}(x,x,y)$, $t(x,x,y) = w_{k}(x,x,y)$, $w_{j-1}(x,x,y) = v_{j}(x,x,y)$, $v_i(x,y,x)$ $ = w_i(x,y,x)$, and $v_i(x,y,y) = w_i(x,y,y)$ are identities of $\vv{V}$, for all $i \in \{1, \dots, k\}$, $j \in \{2, \dots, k\}$. Let $\alg{A} \in \vv{V}$ and let $\delta \in \con(\alg{A})$ with $[\delta,\delta]_{2T} = 0$. We can apply the previous claim to conclude that $v_i^{\alg{A}}(a,b,c)$ $ = w_i^{\alg{A}}(a,b,c)$, for $a,b,c$ in the same $\delta$-class and for all $i \in \{1, \dots, k\}$. Thus, $s^{\alg{A}}(a,a,c) = t^{\alg{A}}(a,a,c)$ with $a \mathrel{\delta} c$.  Moreover, from $t(x,y,x) = s(x,y,x)$ and $s^{\alg{A}}(a,a,c)= t^{\alg{A}}(a,a,c)$ follows that $t^{\alg{A}}(a,b,c) = s^{\alg{A}}(a,b,c)$ for all $a,b,c$ in the same $\delta$-class, applying a symmetric version of the previous claim.  Thus we proved the statement for $(s, t) \in \a \wedge y^1(\alpha,\b,\g)$. If $(s, t) \in \a \wedge y^n(\alpha,\b,\g)$ with $n >1$ the thesis follows applying inductively the previous strategy.
	\end{proof}
	
	Now we are ready to prove the main theorem of the section. The next result provides a bridge between congruence equations that hold in a variety $\vv V$ and congruence equations that hold in the variety $\vv V'$ generated by abelian algebras of the idempotent reduct and thus is interesting for several reasons. Indeed, a systematic production of a Mal'cev condition by congruence equation valid in $\vv V'$ is somehow surprising. This implies also that the commutator equations generated by the algorithm \ref{PixWilAlg} under some mild assumptions produce a weak Mal'cev condition, result that justifies their study.
	%Note that the proof of the next theorem makes deep use of a modified version of what is also called Mal'cev argument (\cite[Lemma $12.1$]{BS.ACIU}). 
	
	\begin{Thm}\label{ThmweakAlg}
	Let $\vv V$ be a variety and suppose that $\a \wedge (\b \circ \g) \leq p(\alpha, \b ,\g)$ is a congruence equation that fails on a $3$-element set with congruences with pair-wise empty intersection, where $p$ is a term for the language $\{\circ, \wedge, \vee\}$. Then the following are equivalent:
	\begin{enumerate}
		\item [(1)] the abelian algebras of the idempotent reduct of $\vv V$ satisfy $\a \wedge (\b \circ \g) \leq p(\alpha, \b ,\g)$;
		\item [(2)] $\vv V$ satisfies the commutator equations  $\Eq_{C}(\a \wedge (\b \circ \g) \leq p^{(k)}(\alpha, \b ,$ $\g))$ for some $k > 2$ (see Algorithm \ref{PixWilAlg});
		%\item [(3)] $\vv V$ satisfies the equations:
		%\begin{equation}\label{commgenEq}
		%	\a \wedge (\b \circ \g) \leq p(\a', \b',\g');
		%\end{equation}
		%where $\alpha' \in  \{\alpha \circ [\beta,\beta] \circ \a, \alpha \circ [\g,\g] \circ \a\}$, $\b' \in  \{\b \circ [\a,\a] \circ \b, \b \circ [\g,\g] \circ \b\}$ and $\g' \in  \{\g \circ %[\a,\a] \circ \g, \g \circ [\b,\b] \circ \g\}$;
		\item [(3)] $\vv V$ satisfies the equations:
		\begin{equation}\label{commgenEq2}
			\a \wedge (\b \circ \g) \leq p(\a', \b',\g');
		\end{equation}
		where $\alpha' = \alpha \circ G(\alpha,\b,\g)\circ \a$, $\b' = \b \circ G(\b,\a,\g)\circ \b$, and $\g' = \alpha \circ G(\g,\a,\b)\circ \a$ with $G(\alpha,\b,\g) \in \{[\beta \wedge (\a \vee \g), \beta \wedge (\a \vee \g)] , [\g \wedge (\a \vee \b),\g \wedge (\a \vee \b)] \}$;

		  %\{\alpha \circ [\beta \wedge (\a \vee \g), \beta \wedge (\a \vee \g)] \circ \a, \alpha \circ [\g \wedge (\a \vee \b),\g \wedge (\a \vee \b)] \circ \a\}$, $\b' \in  \{\b \circ [\a  \wedge (\b \vee \g),\a  \wedge (\b \vee \g)] \circ \b, \b \circ [\g \wedge (\a \vee \b),\g \wedge (\a \vee \b)] \circ \b\}$ and $\g' \in  \{\g \circ [\a  \wedge (\b \vee \g),\a  \wedge (\b \vee \g)] \circ \g, \g \circ [\beta \wedge (\a \vee \g),\beta \wedge (\a \vee \g)] \circ \g\}$;
		\item [(4)] there exist $n \in \NN$ such that $\vv V$ satisfies the equations:
		\begin{equation}\label{herrgenEq}
			\a \wedge (\b \circ \g) \leq  p(\a', \b',\g');
	\end{equation}
		where $\alpha' = \mathrel{\alpha} \circ \mathrel{F(\a,\b, \g)} \circ \mathrel{\alpha} $, $\b' =  \mathrel{\b} \circ \mathrel{F(\b,\a, \g)} \circ  \mathrel{\b} $, and $\g' =  \mathrel{\g} \circ \mathrel{F(\g, \b, \a)} \circ \mathrel{\g}$ with $F(\a,\b, \g) \in \{\b \wedge y^n(\b,\a,\g), \b \wedge y^n(\b,\g,\a), \g \wedge y^n(\g,\b,\a), \g \wedge y^n(\g,\a,\b)\}$.
	\end{enumerate}
	\end{Thm}

	%where $\alpha' = \mathrel{\alpha} \circ \mathrel{\tau_{\a} \wedge \rho_{\a}}\circ \mathrel{\alpha} $, $\b' =  \mathrel{\b} \circ \mathrel{\tau_{\b} \wedge \rho_{\b}} \circ  \mathrel{\b} $, and $\g' =  \mathrel{\g} \circ \mathrel{\tau_{\g} \wedge \rho_{\g}} \circ \mathrel{\g}$ with $(\tau_{\a}, \rho_{\a}) \in \{(\b, \alpha^n), (\b, \g^n) , (\gamma, \b^n), (\gamma, \a^n)\}$, $(\tau_{\b}, \rho_{\b}) \in \{(\a, \b^n), (\a, \g^n) , (\gamma, \b^n), (\gamma, \a^n)\}$, and $(\tau_{\g}, \rho_{\g}) \in $ $\{(\a, $ $\b^n), (\a, \g^n) ,$ $ (\b, $ $\a^n), (\b, \g^n)\}$ where $\alpha^n = y^n(\beta, \alpha, \gamma)$, $\beta^n = y^n(\alpha,$ $ \beta, \gamma)$, and $\gamma^n = x^n(\beta, \alpha, \gamma)$.

	%$[\a,\a]$ with $\bigcup \a \wedge y^n(\a,\b,\g)$ or $\bigcup \a \wedge x^n(\a,\b,\g)$ and the obtained equation holds in $\alg F_{\vv V}(3)$. Permuting the variables we see that we can substitute every occurrence of $[\b,\b]$ with $\bigcup \b \wedge y^n(\b,\a,\g)$ or $\bigcup \b \wedge x^n(\b,\a,\g)$ and of  $[\g,\g]$ with $\bigcup \g \wedge y^n(\g,\b,\a)$ or $\bigcup \g \wedge x^n(\g,\b,\a)$. Since there are finitely many occurrences of $[\a,\a]$, $[\b,\b]$, and $[\g,\g]$ in the equation (\ref{commgenEq}) then there exists an $n \in \NN$ such that the equations (\ref{herrgenEq}) hold in $\alg F_{\vv V}(3)$. Then the equations (\ref{herrgenEq}) hold in $\vv V$ via a standard Mal'cev argument.

	\begin{proof}
		Let us start from $(1) \Rightarrow (2)$.  Let $\alg A \in \vv V$ and let $\theta \in \con(\alg A)$.  Factoring $\alg A$ by $[\theta, \theta]$, we verify that $\Eq(\a \wedge (\b \circ \g) \leq p^{(k)}(\a,\b,\g))$ hold in every $\theta/[\theta, \theta]$-block of the quotient and thus $\Eq_{C}(\a \wedge (\b \circ \g) \leq p^{(k)}(\a,\b,\g))$ hold in $\alg A$, for some $k > 2$. Let $\alg A'$ be the idempotent reduct of $\alg A/[\theta, \theta]$. Since we factored by $[\theta, \theta]$, we can see that every $\theta/[\theta, \theta]$-class of $\alg A/[\theta, \theta]$ is an abelian subalgebra of  $\alg A'$ and thus it has terms satisfying $\Eq(\a \wedge (\b \circ \g) \leq p^{(k)}(\a,\b,\g))$, for some $k > 2$, via the Pixley-Wille algorithm \cite{Pix.LMC, Wil.K}. We can observe that the terms involved in $\Eq(\a \wedge (\b \circ \g) \leq p^{(k)}(\a,\b,\g))$ witness the satisfaction of  $\Eq_{C}(\a \wedge (\b \circ \g) \leq p^{(k)}(\a,\b,\g))$  in $\alg{A}$, since a quotient by $[\theta,\theta]$  has been applied.
		
		For $(2) \Rightarrow (3)$. From the hypothesis we have that there exists $k > 2$ such that $\Eq_{C}(\alpha \wedge (\beta \circ \gamma) \leq p^{(k)}(\a,\b,\g))$ hold in $\vv{V}$. Then let $\mathbf{A} \in \vv{V}$, $a_1 ,a_2 \in A$, and let $\alpha, \b, \g  \in \con(\mathbf{A})$ be such that:
	\begin{equation*}
		(a_1, a_2) \in \a \wedge (\b \circ \g).
	\end{equation*}	
	Hence, there exists $a_3 \in \alg A$ such that $(a_1,a_3) \in \b$ and $(a_3, a_2) \in \gamma$. We want to prove that 
	\begin{equation}
		\label{GoalnewPW}
		(a_1, a_2) \in p^{(k)}(\a', \b',\g').
	\end{equation}
	
	Let $V_{\mathbf{G}(p^{(k)})} = \{z_1,\dots,z_u\}$ be the set of vertices of $\mathbf{G}(p^{(k)})$, let $Tt^3(p^{(k)})$ be defined as in (\ref{eqp}), and let $\{t_1, \dots, t_u\}$ be the set of terms appearing in pairs of $Tt^3(p^{(k)})$. Let $\rho:V_{\mathbf{G}(p^{(k)})} \rightarrow A$ be the assignment defined by:
    \begin{align*}
   &z_1 \mapsto a_1 = t_1^{\alg A}(a_1, a_2, a_3)
    \\&z_2  \mapsto a_2 = t_2^{\alg A}(a_1, a_2, a_3)
    \\&z_i  \mapsto t_i^{\alg{A}}(a_1, a_2, a_3) \text{ for all } 3 \leq i \leq u.
    \end{align*}
    We want to prove that $\rho$ satisfies the hypothesis of $(1)$ in Proposition \ref{PropKK} and this yields (\ref{GoalnewPW}). Let $(z_i, z_j)$ be and edge of $\mathbf{G}(p^{(k)})$ labeled with $X_s$. From the definition of $\Eq_{C}(\a \wedge (\b \circ \g) \leq p^{(k)}(\a,\b,\g))$, we have that $(t_i,t_j) \in Tt_s^3(p^{(k)})$ and
	\begin{equation*}
		t_i(\sigma_s(x_{1}), \sigma_s(x_{2}), \sigma_s(x_{3})) \approx_{C} t_j(\sigma_s(x_{1}), \sigma_s(x_{2}), \sigma_s(x_{3}))
	\end{equation*}     
	is in $\Eq_{C}(\a \wedge (\b \circ \g) \leq p^{(k)}(\a,\b,\g))$, where $\sigma_s$ is a trasversal of $\Eqv(T^3_s(p))$. Suppose that $\alpha$ is the variable substituted to $X_s$ of $p^{(k)}$, i. e. $s =1$. Then, fixing this congruence, we have that $T^3_s(p) = \{(x_1,x_2)\}$ and
	\begin{align*}
		&t_i^{\alg{A}}(a_1, a_2, a_3) \mathrel{\a } t_i^{\alg{A}}(a_{1}, a_{1}, a_{3}) \mathrel{[\b \wedge (\a \vee \g ),\b \wedge (\a \vee \g )]} t_j^{\alg{A}}(a_{1}, a_{1}, a_{3}) 
		\\&t_j^{\alg{A}}(a_{1}, a_{1}, a_{3}) \mathrel{\a}  t_j^{\alg{A}}(a_1, a_2, a_3).
	\end{align*}
	With the same argument we can check that for all $\theta_s \in \{\alpha,\b, \g\}$ congruence substituted to $X_s$ and $\{\tau_s, \delta_s\} = \{\alpha, \b, \g\} \setminus \{\theta_s\}$:
	\begin{align*}\label{commutatoreq}
		&t_i^{\alg{A}}(a_1, a_2, a_3) \mathrel{\theta_s } \overline{t}_i^{\alg{A}}(a_1, a_2, a_3) \\&\overline{t}_i^{\alg{A}}(a_1, a_2, a_3)   \mathrel{[\tau_s \wedge (\theta_s \vee \delta_s),\tau_s \wedge (\theta_s \vee \delta_s)]}\overline{t}_j^{\alg{A}}(a_1, a_2, a_3) \\&\overline{t}_i^{\alg{A}}(a_1, a_2, a_3)  \mathrel{\theta_s}  t_j^{\alg{A}}(a_1, a_2, a_3)
	\end{align*}
	where $\overline{t}_i(x_1,x_2,x_3) = t_i(\sigma_s(x_{1}), \sigma_s(x_{2}), \sigma_s(x_{3}))$, $\overline{t}_j(x_1,x_2,x_3) = t_j(\sigma_s(x_{1})$ $, \sigma_s(x_{2}), \sigma_s(x_{3}))$, and $\sigma_s$ is an appropiate transversal of $\Eqv(T^3_s(p))$. %Let $\rho:Z \rightarrow A$ be the assignment such that $z_1 \mapsto a_1$,$z_2 \mapsto a_2$ and $z_i \mapsto t_i^{\alg{A}}(a_1, a_2, a_3)$ for all $3 \leq i \leq u$. 
    Thus we have that $(t_i^{\alg{A}}(a_1, a_2, a_3), $ $ t_j^{\alg{A}}(a_1, a_2, a_3)) \mathrel{\in} \mathrel{\theta_s} \circ \mathrel{[\tau_s \wedge (\theta_s \vee \delta_s)}, $ $\tau_s \wedge (\theta_s \vee \delta_s)] \circ \mathrel{\theta_s}$ whenever $(z_i,z_j) \in E_{\mathbf{G}(p^{(k)})}$. By $(1)$ of Proposition \ref{PropKK}, we have that $(a_1,a_2) \in  p^{(k)}(\a', \b',\g')$ $ \subseteq p(\a', \b',\g')$.
	
	%$(2) \Rightarrow (4)$ can be proved using the same strategy of $(2) \Rightarrow (3)$. For  $(3) \Rightarrow (2)$ we use the same strategy of \cite{Ker.VWAD}[Lemma $2.7$]. 
	
	For $(3) \Rightarrow (4)$ follows applying Lemma \ref{herringLem}. Indeed let $\alg F_{\vv V}(3)$ the free algebra of $\vv V$ over $3$ generators $\{x,y,z\}$. We see that (\ref{commgenEq2}) hold in $\alg F_{\vv V}(3)$ by the hypothesis. Let $\a = \Cg(\{(x,z)\})$, $\b = \Cg(\{(x,y)\})$, and $\g = \Cg(\{(y,z)\})$.  
	%By \cite{Lip.ACOV}[Theorem $3.4$ $(a) \Rightarrow (b)$] we can substitute every occurrence in  (\ref{commgenEq}) of $[\a,\a]$ with $[\b \vee \g,\a]$ and the obtained equation holds in $\alg F_{\vv V}(3)$. 
	By Lemma \ref{herringLem} we have that $[\a \wedge(\b \vee \g),\a \wedge(\b \vee \g)] \leq \bigcup \a \wedge y^n(\a,\b,\g), \bigcup \a \wedge y^n(\a,\g,\b)$. Thus we can substitute every occurrence in  (\ref{commgenEq2}) of $[\a \wedge(\b \vee \g),\a \wedge(\b \vee \g)]$ with $\bigcup \a \wedge y^n(\a,\b,\g)$ or  $\bigcup \a \wedge y^n(\a,\g,\b)$ and the obtained equation holds in $\alg F_{\vv V}(3)$. Permuting the variables we can substitute every occurrence of $[\beta \wedge (\a \vee \g),\beta \wedge (\a \vee \g)]$ with $\bigcup \b \wedge y^n(\b,\a,\g)$ or $\bigcup \b \wedge y^n(\b,\g,\a)$ and of $[\g \wedge (\a \vee \b),\g \wedge (\a \vee \b)]$ with $\bigcup \g \wedge y^n(\g,\b,\a)$ or $\bigcup \g \wedge y^n(\g,\a,\b)$. Since we have done finitely many substitutions in the equation (\ref{commgenEq2}) then there exists an $n \in \NN$ such that the equations (\ref{herrgenEq}) hold in $\alg F_{\vv V}(3)$. Hence, the equations (\ref{herrgenEq}) hold in $\vv V$ via a standard Mal'cev argument.
	
	$(4) \Rightarrow (1)$.	Let $\vv V'$ be the variety generated by the abelian algebras of the idempotent reduct of $\vv V$. Since the equation (\ref{herrgenEq}) is a congruence equation then, by the Pixley-Wille algorithm \cite{Pix.LMC, Wil.K}, it generates an idempotent Mal'cev condition that holds in $\vv V'$ and we claim that this Mal'cev condition is non-trivial. In fact, $\a \wedge(\b \circ \g)  \leq p(\alpha,\b,\g)$ is non-trivial and it fails on a $3$-element set choosing $\alpha, \b, \gamma$ in a proper way. Furthermore, from the hypothesis, $\alpha, \b, \gamma$, can be chosen with a pairwise empty intersection and thus with $\a \wedge y^n(\a,\b,\g) = 0$ and the same holds permuting $\alpha,\b,\gamma$. Hence, we can see that if $\a \wedge(\b \circ \g)  \leq p(\alpha,\b,\g)$ fails on a $3$-element set, then also the equations (\ref{herrgenEq}) fails choosing the same congruences. %Moreover, as a consequence of \cite{Lip.ACOV}[Theorem $3.4$ (d)],  $\vv V'$ satisfies the equation obtained by (\ref{herrgenEq}) substituting $[\alpha,\a]_{2T}$ in place of $\a \wedge \b^n$, $\a \wedge \g^n$ and their symmetric versions obtained permuting $\alpha,\b,\gamma$, where $[\alpha,\a]_{2T}$ is the two terms commutator defined in \cite{Lip.ACOV} for example.

    Let $\alg{F}_{\vv V'}(\{x,y,z\})$ be the free algebra over three generators of $\vv{V}'$. % generated by the abelian algebras of the idempotent reduct of $\vv{V}$.
    Let $\overline{\a} = \Cg(\{(x,z)\})$, $\overline{\b}  = \Cg(\{(x,y)\})$, and $\overline{\g} = \Cg(\{(y,z)\})$.  Since $\alg{F}_{\vv{V}'}(\{x,y,z\})$ satisfies (\ref{herrgenEq}), by Proposition \ref{PropKK} there exist $l > 2$ and an assignment $\chi$ from $V_{\alg G(p^{(l)})} \rightarrow F_{\vv{V}'}(\{x,y,z\}): x_s \mapsto t_s$ such that for all edges $(x_i, x_j)$ with label $X_k$ of $\alg G(p^{(l)})$, we have $(t_i,t_j) \in \theta_k$, where $\theta_k$ is the congruence substituted to the variable $X_k$ of $p^{(l)}$. Without loss of generality we suppose that $\theta_k = \overline{\a}' = \overline{\a} \circ F(\overline{\a},\overline{\b},\overline{\g}) \circ \overline{\a}$, i.e. $k = 1$. Moreover, suppose that $(t_i, t_j) \in \overline{\a}'$. Then, there exist $s_i, s_j \in \alg{F}_{\vv{V}'}(\{x,y,z\})$ such that $t_i\mathrel{\overline{\a}} s_i \mathrel{F(\overline{\a},\overline{\b},\overline{\g})} s_j \mathrel{\overline{\a}} t_j$ and thus $t_i(x,y,x) \approx s_i(x,y,x)$, $t_j(x,y,x) \approx s_j(x,y,x)$, and $s_i \mathrel{F(\overline{\a},\overline{\b},\overline{\g})} s_j$. We can observe that $s_i$ and $s_j$ satisfy the hypothesis of Lemma \ref{lem:classher} and thus their realization over an algebra $\alg A \in \vv V'$ are equal whenever computed over the same $\delta$-class with $[\delta, \delta]_{2T} = 0$.
 
	Let $\alg{A}$ be an abelian algebra of the idempotent reduct of $\vv{V}$. By \cite[Corollary $4.5$]{KK.TBTC}, $[\phi,\phi] = 0 \Rightarrow [\phi,\phi]_{2T} = 0$ for Taylor varieties since $[\phi, \phi]_{L} = 0 \Rightarrow [\phi,\phi]_{2T} = 0$, where $[\phi,\phi]_{L}$ is the linear commutator and thus $[\phi,\phi]_{2T} = 0$ for all $\phi \in \con(\alg{A})$. 

	%From the hypothesis we have that $\alg{A}$ satisfies $\a \wedge (\b \circ \g) \leq p(\a', \b',\g')$ and 
    Now we prove that $\alg{A}$ satisfies $\a \wedge (\b \circ \g) \leq p(\a, \b,\g)$. %Let $(a,b) \in \a \wedge(\b \circ \g)$, then there exists $c \in A$ such that $(a,b) \in \alpha$, $(a,c) \in \b$, and $(c,b) \in \g$.  We show that  $(a,b) \in p(\a, \b,\g)$. 
    Let $\a, \b, \g \in \con(\alg A)$ and $(a,c) \in \a \wedge (\b \vee \g)$, then $(a,c) \in \a$ and there exists $b \in A$ such that $(a,b) \in \b$ and $(b,c) \in \g$. Let $\psi:V_{\alg G(p^{(l)})} \rightarrow A$ be the assignment obtained by $\chi$ composing it with the realization $t \mapsto t^{\alg A}(a,b,c)$, extending $x_1 \mapsto a$,$x_2 \mapsto c$ such that $x_i \mapsto t_i^{\alg{A}}(a,b,c)$, for all $t_i$ in the codomain of $\chi$. Let us consider $t^{\alg A}_i$ and $t^{\alg A}_j$ such that $t_i(x,y,x) \approx s_i(x,y,x)$, $t_j(x,y,x) \approx s_j(x,y,x)$, and $s_i \mathrel{F(\overline{\a},\overline{\b},\overline{\g})} s_j$. Then we have:
	\begin{equation*}
		t_i^{\alg{A}}(a,b,c) \mathrel{\a} t_i^{\alg{A}}(a,b,a) = s_i^{\alg{A}}(a,b,a) =  s_j^{\alg{A}}(a,b,a) =  t_j^{\alg{A}}(a,b,a) \mathrel{\a} t_j^{\alg{A}}(a,b,c).
	\end{equation*}
	Thus $t^{\alg A}_i(a,b,c) \mathrel{\a} t^{\alg A}_j(a,b,c)$. From this clearly follows that we can cancel the term $F$ appearing in the $p^{(l)}(\a',\b',\g')$, thus reducing it to $p^{(l)}(\a,\b,\g)$. 
    We can prove symmetric relations for $\b$ and $\g$, thus, by Proposition \ref{PropKK}, the assignment $\psi$ witness the satisfaction of $\a \wedge (\b \circ \g) \leq p(\a, \b,\g)$ in $\alg{A}$. Hence, the abelian algebras of idempotent reduct of $\vv{V}$ satisfy $\a \wedge (\b \circ \g) \leq p(\alpha, \b ,\g)$ and $(1)$ holds.
	\end{proof}

We can see that the previous theorem requires to restrict the congruence equation $p \leq q$ that holds in the variety generated by the abelian algebras of the idempotent reduct of a variety to have $p = \a \wedge (\b \circ \g)$. Although quite restrictive, this hypothesis is satisfied by the vast majority of the interesting known Mal'cev conditions characterized by congruence equations (modularity, distributivity, and many others). Moreover, the hypothesis $p = \a \wedge (\b \circ \g)$ can be weaken further as shown in Theorem \ref{TaylorThm} but the difficulty in doing so comes from the constraints given by the commutator equations which run over a set of variables in the same $\theta$-class. This blocks the generalization of the proof of $(2) \Rightarrow (3)$ of Theorem \ref{ThmweakAlg}.
In the more general context of a congruence equation $p \leq q$, where $p$ and $q$ are terms for the language $\{\circ, \wedge, \vee\}$, we can still prove the equivalence of points $(1)$ and $(2)$ of Theorem \ref{ThmweakAlg}. However, without point $(4)$, we cannot find a weak Mal'cev condition that holds in the whole variety $\vv V$ and is equivalent to the property of satisfying a congruence equation in the variety generated by the abelian algebras of the idempotent reduct.

\begin{Thm}\label{ThmweakAlg2}
	Let $\vv V$ be a variety and suppose that $p \leq q$ is a congruence equation, where $p,q$ are terms for the language $\{\circ, \wedge, \vee\}$. Then the following are equivalent:
	\begin{enumerate}
		\item [(1)] the abelian algebras of the idempotent reduct of $\vv V$ satisfy $p \leq q$;
		\item [(2)] $\vv V$ satisfies the commutator equations $\Eq_{C}(p \leq q)$ (see Algorithm \ref{PixWilAlg}).
  \end{enumerate}
\end{Thm}

Let us start from $(1) \Rightarrow (2)$. Let $\alg A \in \vv V$ and let $\theta \in \con(\alg A)$. Similarly to Theorem \ref{ThmweakAlg}, we can factoring $\alg A$ by $[\theta, \theta]$ and verify that $\Eq(p \leq q)$ hold in every $\theta/[\theta, \theta]$-block of the quotient and thus $\Eq_{C}(p \leq q)$ hold in $\alg A$. %Let $\alg A'$ be the idempotent reduct of $\alg A/[\theta, \theta]$. Since we factored by $[\theta, \theta]$ we can see that every $\theta/[\theta, \theta]$-class of $\alg A/[\theta, \theta]$ is an abelian subalgebra of  $\alg A'$ and thus it has terms satisfying $\Eq(p \leq q)$ via the Pixley-Wille algorithm \cite{Pix.LMC, Wil.K}. We can observe that those terms witness the satisfaction of $\Eq_{C}(p \leq q)$  in $\alg{A}$, since a quotient by $[\theta,\theta]$ has been applied.

For $(2) \Rightarrow (1)$ we can observe that, since all the terms given by the Pixley-Wille algorithm are idempotent then, by Remark \ref{idempRem}, the terms witnessing the satisfaction of $\Eq_{C}(p \leq q)$ are idempotent. Let us consider an abelian algebra $\alg A$ of the idempotent reduct of $\vv V$. Since the terms witnessing the satisfaction of $\Eq_{C}(p \leq q)$ are idempotent we have that $\alg A$ satisfies $\Eq_{C}(p \leq q)$. Then the claim follows from the fact that $\alg A$ is abelian.
	
\section*{Conclusions}
In this study, we have explored the relationship between Mal'cev conditions and commutator equations. This interesting connection that generalizes the Mal'cev condition described by the weak difference term, might have many outcomes in the study of the relationship between Mal'cev conditions and commutator properties or TCT-types, as shown in previous results such as \cite{KK.TBTC, KK.TSOC, Lip.ACOV}. 

Our aim with this work was to establish the foundations of a novel type of well-behaved equations that could produce new insights about Mal'cev condition. Further developments could be obtained by applying the procedure outlined in Theorem \ref{ThmweakAlg} to known congruence equations. This approach could reveal whether the resulting commutator equations correspond to existing Mal'cev conditions or if they characterize novel ones. In the latter case, it would be of interest to investigate how they impact the behavior of commutators or the solvability theory in the varieties satisfying them. Examples of prototype results about commutators and solvability theory can be found in \cite[Chapter $6$]{KK.TSOC} for the weak difference term.

Finally, we want to emphasize that this paper primarily focuses on properties satisfied by abelian algebras of the idempotent reduct of a variety. The fact that under mild assumptions those kind of properties generally yield weak Mal'cev conditions shows an unexpected deep connection between a variety and its abelian algebras of the idempotent reduct. This enables the verification of Mal'cev conditions in the more manageable context of abelian algebras. %and, if needed, allows to check Mal'cev conditions in the more manageable setting of abelian algebras.

\section*{Acknowledgement}	

The author thanks Paolo Aglianò, Erhard Aichinger, Sebastian Kreinecker, and Bernardo Rossi for many hours of fruitful discussions.

\bibliographystyle{alpha}

\begin{thebibliography}{HLMS07}
	
	\bibitem{Bul.ADTF}
	Bulatov, A. A.:
	A Dichotomy Theorem for Nonuniform CSPs.
	2017 IEEE 58th Annual Symposium on Foundations of Computer Science (FOCS), 319-330 (2017)
	
	\bibitem{BS.ACIU}
	Burris, S., Sankappanavar, H.P.:
	A Course in Universal Algebra.
	Graduate Texts in Mathematics,
	vol. 78.
	Springer-Verlag, New York-Berlin (1981)
	
	\bibitem{Cze.AMTC}
	Cz\'{e}dli, G.:
	A Mal'cev-type condition for the semi-distributivity of congruence lattices.
	Acta Sci. Math. (Szeged), 
	\textbf{43}(3-4), 267--272 (1981)
	
	\bibitem{CD.HSWW}
	Cz\'{e}dli, G., Day, A.:
	Horn sentences with ($W$) and weak Mal'cev conditions.
	Algebra Universalis,
	\textbf{19}(2), 217--230 (1984)
	
	\bibitem{Day.ACOM}
	Day, A.:
	A characterization of modularity for congruence lattices of algebras.
	Canad. Math. Bull., 
	\textbf{12}, 167--173 (1969)
	
	%\bibitem{Fio.CSOF}
	%Fioravanti, S.:
	%Closed sets of finitary functions between finite fields of coprime order.
	%Algebra Universalis, 
	%\textbf{82}, 52 (2020)
	
	\bibitem{Fio.CSOF2}
	Fioravanti, S.:
	Closed sets of finitary functions between products of finite fields of coprime order.
	Algebra Universalis, 
	\textbf{82}, 61 (2021)
	
	\bibitem{Fio.EOAS}
	Fioravanti, S.:
	Expansions of abelian square-free groups.
	International Journal of Algebra and Computation,
	\textbf{31}, 4  623--638(2021)
	
	\bibitem{Fio.MCCT}
	Fioravanti, S.:
	Mal’cev conditions corresponding to identities for compatible reflexive relations.
	Algebra Universalis, 
	\textbf{82}, 10 (2021)
	
	\bibitem{FM.CTFC}
	Freese, R., McKenzie, R.:
	Commutator theory for congruence modular varieties.
	London Mathematical Society Lecture Note Series,
	
	\bibitem{Jon.AWCL}
	J\'{o}nsson, B.:
	Algebras whose congruence lattices are distributive.
	Math. Scand., 
	\textbf{21}, 110--121 (1967)
	
	\bibitem{Ker.VWAD}
	Kearnes, K.A.:
	Varieties with a difference term.
	J. Algebra,
	177, pp. 926--960 (1995)
	
	\bibitem{Ker.RMDO}
	Kearnes, K.A.:
	Relative Maltsev definability of some commutator properties.
	https://arxiv.org/pdf/2202.10009.pdf
	
	
	\bibitem{KK.TBTC}
	Kearnes, K.A., Szendrei,  Å.:
	The relationship between two commutators.
	International Journal of Algebra and Computations, 
	\textbf{8}(4), 497--531 (1998)
	
	\bibitem{KK.TSOC}
	Kearnes, K.A., Kiss, E.W.:
	The shape of congruence lattices.
	Mem. Amer. Math. Soc., 
	\textbf{222}(1046) (2013)
	
	\bibitem{Lip.ACOV}
	Lipparini, P.:
	A characterization of varieties with a difference term.
	Canadian Mathematical Bulletin, 
	\textbf{29}(3), 308--315 (1996)
	
	\bibitem{Lip.FCIT}
	Lipparini, P.:
	From congruence identities to tolerance identities.
	Acta Sci. Math. (Szeged), 
	\textbf{73}(1-2), 31--51 (2007)
	
	\bibitem{Mal.OTGT}
	Mal'cev, A.I.:
	On the general theory of algebraic systems.
	Mat. Sb. N.S., 
	\textbf{35}(77),3--20 (1954) (Russian)
	
	\bibitem{Ols.TWNI}
	Ol\v{s}\'{a}k, M.,
	The weakest nontrivial idempotent equations.
	Bulletin of the London Mathematical Society,
	\textbf{49}(6),  1028--1047 (2017)
	
	\bibitem{Pix.DAPO}
	Pixley, A.F.:
	Distributivity and permutability of congruence relations in
	equational classes of algebras.
	Proc. Amer. Math. Soc., 
	\textbf{14}, 105--109 (1963)
	
	\bibitem{Pix.LMC}
	Pixley, A.F.:
	Local Mal'cev conditions.
	Canad. Math. Bull., 
	\textbf{15}, 559--568 (1972)

    \bibitem{Spa.OTNC}
	Sparks, A.:
	On the number of clonoids.
	Algebra Universalis
	\textbf{80}(4), 53 (2019)
	
	\bibitem{Taylor1973}
	Taylor, W:
	Characterizing {M}al'cev conditions, 
	Algebra Universalis {\bf
		3},  351--397 (1973)
	
	\bibitem{Tay.VOHL}
	Taylor, W.:
	Varieties Obeying Homotopy Laws.
	Canadian Journal of Mathematics, 
	\textbf{3}, 498--527, Cambridge University Press (1977)

	
	\bibitem{Wil.K}
	Wille, R.:
	Kongruenzklassengeometrien.
	Lecture Notes in Mathematics, vol. 113. Springer-Verlag, Berlin-New York (1970)
	
	\bibitem{Zhuk.APOT}
	Zhuk, D.:
	A Proof of the {CSP} Dichotomy Conjecture.
	J. {ACM}, 
	\textbf{67}(5), 30:1--30:78 (2020)
\end{thebibliography}
\end{document}